\documentclass[12pt]{article}%
\usepackage{graphicx}
\usepackage[intlimits]{amsmath}
\usepackage{latexsym}
\usepackage{amsfonts}
\usepackage{amssymb}%
 \makeatletter
\newtheorem{theorem}{Theorem}

\newtheorem{fact}[theorem]{Fact}
\newtheorem{lemma}[theorem]{Lemma}

\newenvironment{proof}[1][Proof]{\noindent{\textbf {#1}  }}  {\hfill$\Box$\bigskip}

\begin{document}

\title{Large joints in graphs}
\author{B\'{e}la Bollob\'{a}s \thanks{Department of Pure Mathematics and Mathematical
Statistics, University of Cambridge, Cambridge CB3 0WB, UK and}
\thanks{Department of Mathematical Sciences, University of Memphis, Memphis TN
38152, USA} \thanks{Research supported in part by NSF grants DMS-0505550,
CNS-0721983 and CCF-0728928, and ARO grant W911NF-06-1-0076} \ and Vladimir
Nikiforov \thanks{Research supported by NSF Grant \# DMS-0906634}
\footnotemark[2]}
\maketitle

\begin{abstract}
We show that if $r\geq s\geq2,$ $n>r^{8},$ and $G$ is a graph of order $n$
containing as many $r$-cliques as the $r$-partite Tur\'{a}n graph of order
$n,$ then $G$ has more than at $n^{r-1}/\left(  4r\right)  ^{r+6}$ cliques
sharing a common edge unless $G$ is isomorphic to the the $r$-partite
Tur\'{a}n graph of order $n$. This structural result generalizes a previous
result that has been useful in extremal combinatorics.\medskip\ 

\textbf{Keywords: }\textit{joint; jointsize; clique; number of cliques;
Tur\'{a}n graph }

\end{abstract}

\section*{Introduction}

In notation we follow \cite{Bol98}; in particular, $T_{r}\left(  n\right)  $
denotes the $r$-partite Tur\'{a}n graph of order $n$ and $t_{r}\left(
n\right)  $ denotes the number of its edges. Also, an $r$-joint of size $t$ is
a collection of $t$ distinct $r$-cliques sharing an edge. (Note that two
$r$-cliques of an $r$-joint may share $r-1$ vertices.) We write $\mathrm{js}%
_{r}\left(  G\right)  $ for the maximum size of an $r$-joint in a graph $G$;
in particular, if $2\leq r\leq n$ and $r$ divides $n$ then $\mathrm{js}%
_{r}(K_{n})=\binom{n-2}{r-2}$ and $\mathrm{js}_{r}\left(  T_{r}(n)\right)
=\left(  \frac{n}{r}\right)  ^{r-2}$.

In \cite{BoNi04} we improved a result of Erd\H{o}s \cite{Erd69} to the
following assertion.\medskip

\emph{Let $r\geq2$, $n>r^{8}$, and let $G$ be a graph of order $n$ and size at
least $t_{r}(n)$. Then
\begin{equation}
\mathrm{js}_{r+1}\left(  G\right)  >\frac{n^{r-1}}{r^{r+5}}\label{ourb0}%
\end{equation}
unless $G=T_{r}\left(  n\right)  $.}\medskip

Joints have a long history in graph theory. The study of $\mathrm{js}%
_{3}\left(  G\right)  $, also known as the \emph{booksize} of $G$, was
initiated by Erd\H{o}s in \cite{Erd62b} and subsequently generalized in
\cite{Erd68} and \cite{Erd69}; it seems that he foresaw the importance of
joints when he restated his general results in 1995, in \cite{EFGG95}. A
quintessential result concerning joints is the \textquotedblleft triangle
removal lemma\textquotedblright\ of Ruzsa and Szemer\'{e}di \cite{RuSz75},
which can be stated as a lower bound on $\mathrm{js}_{3}\left(  G\right)  $
when $G$ is a graph of a particular kind. Erd\H{o}s's challenge was taken up
in Ramsey graph theory, see, e.g., \cite{NiRo09} and its references. Some
recent applications are given in \cite{Nik09b, Nik09c}.

\vspace{5pt} Our aim in this note is to prove an analogue of inequality
\eqref{ourb0} in the case when $G$ has a fair number of $r$-cliques, rather
than edges. More precisely, letting $k_{r}\left(  G\right)  $ stand for the
number of $r$-cliques of a graph $G$, we shall prove the following theorem.

\begin{theorem}
\label{th1} Let $r\geq s\geq2,$ $n>r^{8},$ and let $G$ be a graph of order
$n$, with $k_{s}(G) \geq k_{s}\left(  T_{r}(n) \right)  $. Then
\begin{equation}
\label{new}\mathrm{js}_{r+1}\left(  G\right)  >\frac{n^{r-1}}{\left(
4r\right)  ^{r+6}}%
\end{equation}
unless $G=T_{r}\left(  n\right)  $.
\end{theorem}

Inequality (\ref{new}) is far from the best possible; in particular, for $s=2$
inequality (\ref{ourb0}) is significantly better. However, in most
applications, the exact values of the coefficients to $n^{r-1}$ in
(\ref{ourb0}) and (\ref{new}) are irrelevant, except for the convenience.
Moreover, these inequalities cannot be improved too much, as shown by the
graph $G$ obtained by adding an edge to $T_{r}(n)$: if $n$ is a multiple of
$r$ then $k_{s}(G)\geq k_{s}(T_{r}(n))$ and $\mathrm{js}_{r+1}(G)=\left(
\frac{n}{r}\right)  ^{r-1}$.

\vspace{5pt} We very much hope that Theorem \ref{th1} will be one of many new
generalizations of classical extremal results in graph theory to be proved in
the near future.

\section*{Preliminary results}

In this section we shall collect the results we shall use in our proof of
Theorem~\ref{th1}. The first two, stated as `Facts', are from earlier papers,
but the two lemmas following them seem to be new. We shall also need two
simple inequalities about the Tur\'{a}n graph $T_{r}(n)$. The required proofs
of the results below will be given in the next section.

We start with an inequality stated by Moon and Moser in \cite{MoMo62}; it
seems that Khad\v{z}iivanov and Nikiforov \cite{KhNi78} were the first to
publish a complete proof of this (see also \cite{Lov79}, Problem 11.8).

\begin{fact}
Let $1\leq s<t<n$, and let $G$ be a graph of order $n$ containing at least one
$t$-clique. Then
\begin{equation}
\frac{\left(  t+1\right)  k_{t+1}\left(  G\right)  }{tk_{t}\left(  G\right)
}-\frac{n}{t}\geq\frac{\left(  s+1\right)  k_{s+1}\left(  G\right)  }%
{sk_{s}\left(  G\right)  }-\frac{n}{s}. \label{MoMo}%
\end{equation}

\end{fact}

\vspace{10pt}

The second fact we need is a stability theorem, stated as Theorem 9 in
\cite{BoNi04}.

\begin{fact}
\label{stabj} Let
\[
r\geq2,\text{ \ \ }n>r^{8} \ \text{and } \ 0<\beta<r^{-8}/16;
\]
furthermore, let $G$ be a graph of order $n$ and size
\[
e\left(  G\right)  >\left(  \frac{r-1}{2r}-\beta\right)  n^{2}.
\]
Then either
\begin{equation}
\mathrm{js}_{r+1}\left(  G\right)  >\frac{n^{r-1}}{r^{r+6}}, \label{minjs}%
\end{equation}
or $G$ contains an induced $r$-partite subgraph $G_{0}$ of order at least
$\left(  1-2\sqrt{\beta}\right)  n$ with minimum degree
\begin{equation}
\delta\left(  G_{0}\right)  >\left(  1-\frac{1}{r}-4\sqrt{\beta}\right)  n.
\label{mindg}%
\end{equation}

\end{fact}

Let us turn to the two technical lemmas, which seem to be new. The first one
is somewhat paradoxical: informally it says that if a graph $G$ contains few
$\left(  r+1\right)  $-cliques, then the ratio $k_{2}\left(  G\right)
/k_{r}\left(  G\right)  $ is large.

\begin{lemma}
\label{le1}Let $\alpha\geq0$ and $G$ be a graph of order $n.$ If
\[
k_{r+1}\left(  G\right)  <\frac{\alpha r^{2}}{r+1}\left(  \frac{n}{r}\right)
^{r+1},
\]
then
\[
k_{2}\left(  G\right)  >\frac{rk_{r}\left(  G\right)  }{2n^{r-2}}\prod
_{s=2}^{r-1}\left(  \frac{r-s}{rs}+\alpha\right)  ^{-1}.
\]

\end{lemma}

\begin{lemma}
\label{le2} Let $\alpha>0,$ $2\leq s\leq r\leq n,$ and let $G$ be a graph of
order $n.$ If $k_{r}\left(  G\right)  \geq k_{r}\left(  T_{r}\left(  n\right)
\right)  ,$ then either
\[
\mathrm{js}_{r+1}\left(  G\right)  >\alpha r\left(  \frac{n}{r}\right)  ^{r-1}%
\]
or
\[
k_{2}\left(  G\right)  >\left(  \frac{r-1}{2r}-\frac{r^{3}\alpha}{2}%
-\frac{r^{3}}{16n^{2}}\right)  n^{2}.
\]

\end{lemma}

Finally, the following two inequalities about Tur\'an graphs are easily checked.

\begin{fact}
For every $2\leq r\leq n,$
\begin{align}
k_{2}\left(  T_{r}\left(  n\right)  \right)   &  \geq\frac{r-1}{2r}n^{2}%
-\frac{r}{8}.\label{tur2}\\
k_{r}\left(  T_{r}\left(  n\right)  \right)   &  \geq\left(  \frac{n}%
{r}\right)  ^{r}-\frac{r^{2}}{16}\left(  \frac{n}{r}\right)  ^{r-2}.
\label{turr}%
\end{align}

\end{fact}

\section*{Proofs}

In this section we shall prove Lemmas \ref{le1} and \ref{le2}, and Theorem
\ref{th1}.

\vspace{5pt}

\begin{proof}
[\textbf{Proof of Lemma \ref{le1}}.]We have
\[
\frac{\left(  r+1\right)  k_{r+1}\left(  G\right)  }{rk_{r}\left(  G\right)
}<\alpha r\left(  \frac{n}{r}\right)  ^{r+1}\binom{n}{r}^{-1}\leq\alpha
r\left(  \frac{n}{r}\right)  ^{r+1}\left(  \frac{n}{r}\right)  ^{-r}=\alpha
n.
\]
Now, for every $s=2,\ldots,r-1,$ inequality (\ref{MoMo}) gives
\[
\frac{\left(  s+1\right)  k_{s+1}\left(  G\right)  }{sk_{s}\left(  G\right)
}-\frac{n}{s}\leq\frac{\left(  r+1\right)  k_{r+1}\left(  G\right)  }%
{rk_{r}\left(  G\right)  }-\frac{n}{r}\leq\alpha n-\frac{n}{r},
\]
and so,%
\[
\frac{\left(  s+1\right)  k_{s+1}\left(  G\right)  }{sk_{s}\left(  G\right)
}\leq\left(  \frac{r-s}{sr}+\alpha\right)  n.
\]
Multiplying these inequalities for $s=2,\ldots,r-1,$ we obtain%
\[
2k_{2}\left(  G\right)  n^{r-2}\prod_{s=2}^{r-1}\left(  \frac{r-s}{rs}%
+\alpha\right)  \geq rk_{r}\left(  G\right)  ,
\]
and the desired inequality follows.
\end{proof}

\begin{proof}
[\textbf{Proof of Lemma \ref{le2}}.]Assume that $\mathrm{js}_{r+1}\left(
G\right)  \leq\alpha r\left(  n/r\right)  ^{r-1}.$ Then
\[
\binom{r+1}{2}k_{r+1}\left(  G\right)  \leq\mathrm{js}_{r+1}\left(  G\right)
k_{2}\left(  G\right)  <\alpha r\left(  \frac{n}{r}\right)  ^{r-1}\frac{n^{2}%
}{2}%
\]
and so,%
\[
k_{r+1}\left(  G\right)  \leq\alpha\frac{r^{2}}{r+1}\left(  \frac{n}%
{r}\right)  ^{r+1}.
\]
Now Lemma \ref{le1} and inequality (\ref{turr}) give%
\begin{align*}
k_{2}\left(  G\right)   &  >\frac{rk_{r}\left(  G\right)  }{2n^{r-2}}%
\prod_{i=2}^{r-1}\left(  \frac{r-i}{ri}+\alpha\right)  ^{-1}\\
&  >r\left(  \frac{1}{r}\right)  ^{r-2}\left(  \left(  \frac{n}{r}\right)
^{2}-\frac{r^{2}}{16}\right)  \prod_{i=2}^{r-1}\left(  \frac{r-i}{ri}%
+\alpha\right)  ^{-1}.
\end{align*}
Furthermore, note that%
\begin{align*}
\prod_{i=2}^{r-1}\left(  \frac{r-i}{ri}+\alpha\right)   &  =\prod_{i=2}%
^{r-1}\left(  1+\frac{ri}{r-i}\alpha\right)  \prod_{i=2}^{r-1}\left(
\frac{r-i}{ri}\right)  \leq\left(  1+r\left(  r-1\right)  \alpha\right)
^{r-2}\prod_{i=2}^{r-1}\left(  \frac{r-i}{ri}\right) \\
&  =\left(  1+r\left(  r-1\right)  \alpha\right)  ^{r-2}\left(  \frac{1}%
{r}\right)  ^{r-2}\frac{r-2}{2}\cdot\frac{r-3}{3}\cdot\cdots\cdot\frac{2}%
{r-2}\cdot\frac{1}{r-1}\\
&  =\left(  \frac{1}{r}\right)  ^{r-2}\frac{1}{r-1}\left(  1+r\left(
r-1\right)  \alpha\right)  ^{r-2}.
\end{align*}
Hence, by (\ref{tur2}), we see that
\begin{align*}
k_{2}\left(  G\right)   &  >\binom{r}{2}\frac{1}{\left(  1+r\left(
r-1\right)  \alpha\right)  ^{r-2}}\left(  \left(  \frac{n}{r}\right)
^{2}-\frac{r^{2}}{16}\right) \\
&  >\left(  \frac{r-1}{2r}\right)  \left(  1-r^{2}\alpha\right)  ^{r-2}\left(
1-\frac{r^{4}}{16n^{2}}\right)  n^{2}\\
&  >\left(  \frac{r-1}{2r}-\frac{r^{3}\alpha}{2}-\frac{r^{4}}{16n^{2}}\right)
n^{2}.
\end{align*}
as claimed.
\end{proof}

\vspace{10pt} After all this preparation, we are ready to prove our main result.

\vspace{5pt}

\begin{proof}
[\textbf{Proof of Theorem \ref{th1}}.]As shown in \cite{Bol76} (see also
\cite{Bol78}, p.??), if $k_{s}(G)\geq k_{s}\left(  T_{r}(n)\right)  ,$ then
$k_{r}(G)\geq k_{r}\left(  T_{r}(n)\right)  $. Consequently, we may assume
that $s=r.$ Also, assume for a contradiction that
\begin{equation}
\mathrm{js}_{r+1}\left(  G\right)  \leq\frac{n^{r-1}}{\left(  4r\right)
^{r+6}}. \label{assum}%
\end{equation}

First, setting $\alpha=4^{-r-6}r^{-7},$ Lemma \ref{le2} implies that
\[
e\left(  G\right)  >\left(  \frac{r-1}{2r}-\frac{r^{3}\alpha}{2}-\frac{r^{4}%
}{16n^{2}}\right)  n^{2}>\left(  \frac{r-1}{2r}-\frac{1}{4r^{12}}\right)
n^{2}.
\]
Now, recalling that $n>r^{8}$ and setting $\beta=r^{-12}/4,$ by
Fact~\ref{stabj} we find that $G$ contains an induced $r$-partite subgraph
$G_{0}$ with $\left\vert G_{0}\right\vert \geq\left(  1-r^{-6}\right)  n$ and
minimum degree $\delta\left(  G_{0}\right)  >\left(  1-1/r-2r^{-6}\right)  n.$

Let $V_{1},\ldots,V_{r}$ be the vertex classes of $G_{0},$ set $V_{0}=V\left(
G\right)  \backslash V\left(  G_{0}\right)  ,$ and let $U$ be the set of
vertices in $V_{0}$ joined to a vertex of each $V_{1},\ldots,V_{r}.$ Set for
short $\varepsilon=2r^{-6}$ and $\delta=\delta\left(  G_{0}\right)  .$ It
turns out that none of the vertex classes is significantly larger than $n/r.$
Indeed, for every $i\in\left[  r\right]  ,$ we see that%
\begin{equation}
\left\vert V_{i}\right\vert \leq\left\vert G_{0}\right\vert -\delta\leq
n-\left(  1-1/r-\varepsilon\right)  n=\left(  \frac{1}{r}+\varepsilon\right)
n. \label{ubv}%
\end{equation}
Before giving further details, we shall outline the remaining steps of our
proof in three formal claims.

\textbf{Claim 1.} \emph{For every }$u\in U,$\emph{ there exist two distinct
elements }$i,j\in\left[  r\right]  $\emph{ such that}
\[
\left\vert \Gamma(u) \cap V_{i}\right\vert <\frac{n}{3r} \ \ \ \text{and}
\ \ \ \left\vert \Gamma(u) \cap V_{j}\right\vert <\frac{n}{3r}.
\]

\textbf{Claim 2.}\emph{ Every vertex }$u\in U$\emph{ belongs to at most
}$0.91\left(  n/r\right)  ^{r-1}$\emph{ distinct }$r$\emph{-cliques of }%
$G$\emph{.}$\medskip$

\textbf{Claim 3.}\emph{ If }$U$\emph{ is non-empty, then }$k_{r}\left(
G\right)  <k_{r}\left(  T_{r}\left(  n\right)  \right)  $.$\medskip$

The last Claim gives us a contradiction if $U\neq\emptyset$. However, if $U$
is empty, the graph $G$ is $r$-partite and $k_{r}(G)\leq k_{r}\left(
T_{r}(n)\right)  $, with equality if and only if $G=T_{r}\left(  n\right)  $.
Hence, to complete our proof of Theorem~\ref{th1}, all that remains is to
prove these claims.$\medskip$

\textbf{Proof of Claim 1. }Assume for a contradiction that there is $u\in U$
such that
\[
\left\vert V_{i}\cap\Gamma\left(  u\right)  \right\vert \geq\frac{1}{3r}n
\]
for all but at most one $i\in\left[  r\right]  $; if there is such an $i$, we
may assume that $i=1$. Choose $v_{1}\in V_{1}\cap\Gamma\left(  u\right)  ;$ we
shall prove that the edge $uv_{1}$ is contained in at least $\left(
n/4r\right)  ^{r-1}$ distinct $\left(  r+1\right)  $-cliques. This will give
$\mathrm{js}_{r+1}\left(  G\right)  \geq\left(  n/4r\right)  ^{r-1},$
contradicting the assumption (\ref{assum}).

Let $2\leq s\leq r-1$ and choose any $s-1$ vertices $v_{i}\in V_{i},$
$i\in\left[  2..r\right]  .$ Letting%
\[
b=\left\vert V_{s+1}\cap\Gamma\left(  u\right)  \cap\left(  \cap_{i=1}%
^{s}\Gamma\left(  v_{i}\right)  \right)  \right\vert ,
\]
we shall prove that $b>n/\left(  4r\right)  .$ Indeed, for every $i\in\left[
2..s\right]  ,$ note that%
\begin{align*}
\left\vert V_{s+1}\cap\Gamma\left(  v_{i}\right)  \right\vert  &  =\left\vert
V_{s+1}\right\vert +\left\vert \Gamma\left(  v_{i}\right)  \right\vert
-\left\vert V_{s+1}\cup\Gamma\left(  v_{i}\right)  \right\vert \geq\left\vert
V_{s+1}\right\vert +\delta-\left(  n^{\prime}-\left\vert V_{i}\right\vert
\right)  \\
&  =\left\vert V_{s+1}\right\vert +\delta-n+\left\vert V_{i}\right\vert .
\end{align*}
Now, we find that
\begin{align*}
b &  =\left\vert V_{s+1}\cap\Gamma\left(  u\right)  \cap\left(  \cap_{i=1}%
^{s}\Gamma\left(  v_{i}\right)  \right)  \right\vert \\
&  \geq\left\vert V_{s+1}\cap\Gamma\left(  u\right)  \right\vert +\left\vert
\cap_{i=1}^{s}\left(  V_{s+1}\cap\Gamma\left(  v_{i}\right)  \right)
\right\vert -\left\vert V_{s+1}\right\vert \\
&  \geq\frac{1}{3r}n+\left(  \sum_{i=2}^{s}\left\vert V_{s+1}\cap\Gamma\left(
v_{i}\right)  \right\vert -\left(  s-1\right)  \left\vert V_{s+1}\right\vert
\right)  -\left\vert V_{s+1}\right\vert \\
&  \geq\frac{1}{3r}n+\sum_{i=2}^{s}\left(  \left\vert V_{s+1}\right\vert
+\delta-n+\left\vert V_{i}\right\vert \right)  -s\left\vert V_{s+1}\right\vert
>\frac{1}{2r}n+\sum_{i=2}^{s}\left(  \delta+\left\vert V_{i}\right\vert
-n\right)  \\
&  >\frac{1}{3r}n+\sum_{i=1}^{r}\left(  \delta+\left\vert V_{i}\right\vert
-n\right)  =\frac{1}{2r}n+r\delta+n^{\prime}-rn\\
&  >\frac{1}{3r}n+\left(  r-1-r\varepsilon\right)  n+\left(  1-\varepsilon
\right)  n-rn>`\left(  \frac{1}{3r}-\left(  r+1\right)  \varepsilon\right)
n\\
&  >\frac{1}{4r}n.
\end{align*}

To bound the number of cliques containing $uv_{1}$, for $s=2,\ldots,r,$ choose
a vertex $v_{s}$ such that
\[
v_{s}\in V_{s}\cap\Gamma\left(  u\right)  \cap\left(  \cap_{i=1}^{s-1}%
\Gamma\left(  v_{i}\right)  \right)  .
\]
Clearly for every choice of $v_{2},\ldots,v_{r},$ the set $\{u,v_{1},v_{2},
\ldots, v_{r}\}$ induces an $\left(  r+1\right)  $-clique. Since for every
$s=2,\ldots,r,$ the vertex $v_{s}$ can be chosen in at least $n/\left(
4r\right)  $ ways, there are at least $\left(  n/4r\right)  ^{r-1}$ distinct
$\left(  r+1\right)  $-cliques containing the edge $uv_{1},$ completing the
proof of Claim 1. \hfill{$\square$}

\vspace{5pt} \textbf{Proof of Claim 2.}\emph{ }Fix a vertex $u\in U$ and let
$K$ be the set of all $r$-cliques containing $u.$ By Claim 1, we can assume
that
\[
\left\vert \Gamma\left(  u\right)  \cap V_{1}\right\vert <\frac{n}%
{3r}\ \ \ \text{and \ \ }\left\vert \Gamma\left(  u\right)  \cap
V_{2}\right\vert <\frac{n}{3r}.
\]
Write $K_{s}$ for the set of $\left(  r-1\right)  $-cliques in $K$
intersecting $V\left(  G_{0}\right)  $ in exactly $s$ vertices and note that%
\[
\left\vert K\right\vert =\left\vert K_{0}\right\vert +\left\vert
K_{1}\right\vert +\cdots+\left\vert K_{r-1}\right\vert .
\]
Since $G_{0}$ is $r$-partite and each vertex class satisfies (\ref{ubv}), for
every $s=1,\ldots,r-1,$
\[
k_{s}\left(  G_{0}\right)  \leq\binom{r}{s}\left(  \frac{1}{r}+\varepsilon
\right)  ^{s}n^{s}.
\]
On the other hand, for $s=1,\ldots,r-1$ there are at most $\binom{\varepsilon
n}{s}$ $s$-cliques entirely outside $G_{0}.$ Thus, for every $s=1,\ldots,r-2,$
we have%
\[
\left\vert K_{s}\right\vert <\binom{\varepsilon n}{r-1-s}\binom{r}{s}\left(
\frac{1}{r}+\varepsilon\right)  ^{s}n^{s}.
\]
It is easy to check that the right-hand side of this inequality increases with
$s,$ and so%
\begin{align}
\left\vert K_{1}\right\vert +\cdots+\left\vert K_{r-1}\right\vert  &  <\left(
r-1\right)  \varepsilon n\binom{r}{r-2}\left(  \frac{1}{r}+\varepsilon\right)
^{r-2}n^{r-2}\label{cli}\\
&  <\frac{r^{3}}{2}\varepsilon\left(  \frac{1}{r}+\varepsilon\right)
^{r-2}n^{r-1}<\frac{1}{r^{3}}\left(  \frac{1}{r}+\varepsilon\right)
^{r-2}n^{r-1}\nonumber\\
&  <\frac{1}{r^{2}}\left(  \frac{1}{r}+\varepsilon\right)  ^{r-1}n^{r-1}%
\leq\frac{1}{9}\left(  \frac{1}{r}+\varepsilon\right)  ^{r-1}n^{r-1}.\nonumber
\end{align}

Looking closely at $K_{0},$ it turns out that $K_{0}$ is the union of the
following three disjoint sets:
\begin{align*}
K_{0}^{\prime}  &  =\left\{  R:R\in K_{0},\text{ }R\cap V_{1}\neq
\emptyset,\text{ }R\cap V_{2}=\emptyset\right\}  ,\\
K_{0}^{\prime\prime}  &  =\left\{  R:R\in K_{0},\text{ }R\cap V_{2}%
\neq\emptyset,\text{ }R\cap V_{1}=\emptyset\right\}  ,\\
K_{0}^{\prime\prime\prime}  &  =\left\{  R:R\in K_{0},\text{ }R\cap V_{1}%
\neq\emptyset,\text{ }R\cap V_{2}\neq\emptyset\right\}  .
\end{align*}
Thus,
\begin{align*}
\left\vert K_{0}\right\vert  &  =\left\vert K_{0}^{\prime}\right\vert
+\left\vert K_{0}^{\prime\prime}\right\vert +\left\vert K_{0}^{\prime
\prime\prime}\right\vert \leq2\frac{1}{3r}n\left(  \frac{1}{r}+\varepsilon
\right)  ^{r-3}n^{r-2}+\frac{1}{9r^{2}}n^{2}\left(  \frac{1}{r}+\varepsilon
\right)  ^{r-3}n^{r-3}\\
&  =\frac{2}{3}\left(  \frac{1}{r}+\varepsilon\right)  ^{r-1}n^{r-1}+\frac
{1}{9}\left(  \frac{1}{r}+\varepsilon\right)  ^{r-1}n^{r-1}=\frac{7}{9}\left(
\frac{1}{r}+\varepsilon\right)  ^{r-1}n^{r-1}.
\end{align*}
Hence, in view of (\ref{cli}),
\[
\left\vert K\right\vert <\frac{8}{9}\left(  \frac{1}{r}+\frac{2}{r^{6}%
}\right)  ^{r-1}n^{r-1}<\frac{8}{9}\left(  1+\frac{2}{r^{5}}\right)
^{r-1}\left(  \frac{n}{r}\right)  ^{r-1}<0.91\left(  \frac{n}{r}\right)
^{r-1},
\]
completing the proof of Claim 2. \hfill{$\square$}

\vspace{5pt} \textbf{Proof of Claim 3.} First note that%
\begin{align*}
k_{r}\left(  T_{r}\left(  n\right)  \right)  -k_{r}\left(  T_{r}\left(
n-1\right)  \right)   &  \geq\left(  \frac{n}{r}-1\right)  ^{r-1}>\left(
\frac{n}{r}\right)  ^{r-1}\left(  1-\frac{r\left(  r-1\right)  }{n}\right) \\
&  >\left(  1-\frac{r\left(  r-1\right)  }{r^{8}}\right)  \left(  \frac{n}%
{r}\right)  ^{r-1}\\
&  >\left(  \frac{n}{r}\right)  ^{r-1}\left(  1-\frac{2}{3^{7}}\right)
>0.99\left(  \frac{n}{r}\right)  ^{r-1}.
\end{align*}
In particular, this implies that
\[
k_{r}\left(  T_{r}\left(  n\right)  \right)  -k_{r}\left(  T_{r}\left(
n-\left\vert U\right\vert \right)  \right)  >0.99\left\vert U\right\vert
\left(  \frac{n-\left\vert U\right\vert }{r}\right)  ^{r-1}.
\]
According to Claim 2, by removing the set $U$ we destroy at most
\[
0.91\left\vert U\right\vert \left(  \frac{n}{r}\right)  ^{r-1}%
\]
$r$-cliques. But the graph induced by $V\left(  G\right)  \backslash U$ is
$r$-partite and so, according to Zykov's theorem, \cite{Zyk49}, it has at most
$k_{r}\left(  T_{r}\left(  n-\left\vert U\right\vert \right)  \right)  $
$r$-cliques. Thus,%
\[
k_{r}\left(  T_{r}\left(  n-\left\vert U\right\vert \right)  \right)
+0.91\left\vert U\right\vert \left(  \frac{n}{r}\right)  ^{r-1}\geq
k_{r}\left(  G\right)  \geq k_{r}\left(  T_{r}\left(  n\right)  \right)  ,
\]
implying in turn that
\[
0.91\left(  \frac{n}{r}\right)  ^{r-1}>0.99\left(  \frac{n-\left\vert
U\right\vert }{r}\right)  ^{r-1},
\]
and so%
\[
\frac{0.91}{0.99}>\left(  1-\frac{\left\vert U\right\vert }{n}\right)
^{r-1}>\left(  1-\frac{1}{r^{6}}\right)  ^{r-1}>\left(  1-\frac{1}{3^{6}%
}\right)  ^{2}.
\]
This contradiction completes the proof of Claim 3 and Theorem \ref{th1}.
\end{proof}

\vspace{5pt} It would be good to determine the best constant in
Theorem~\ref{new}, the maximal $c$ such that if $2\leq s\leq r$ are fixed,
$n\rightarrow\infty$, and $G$ is a graph of order $n$ with $k_{s}(G)\geq
k_{s}(T_{r}(n))$ then $\mathrm{js}_{r+1}(G)\geq\left(  c+o(1)\right)  n^{r-1}$
unless $G=T_{r}(n)$. For $s=r=2,$ it is known that the best constant is $1/6$
see \cite{BoNi05} and the references therein. For larger values of $r$, we do
not expect this task to be easy.

\end{document}